\documentclass[preprint,12pt,authoryear]{elsarticle}

\usepackage[letterpaper,margin=1in]{geometry}
\usepackage{amsmath,amssymb,amsthm,mathtools}
\usepackage{mathrsfs}
\biboptions{authoryear,round}
\usepackage{url}
\usepackage{adjustbox}
\usepackage{array,booktabs,longtable}
\usepackage{fvextra}
\DefineVerbatimEnvironment{WLExample}{Verbatim}{fontsize=\footnotesize,breaklines=true,breakanywhere=true,breaksymbolleft={},breaksymbolright={}}
\usepackage[colorlinks=true,linkcolor=blue,citecolor=blue,urlcolor=blue]{hyperref}
\hypersetup{
  pdftitle={Certified Arbitrary-Precision Evaluation of a Family of Generalized Multiple Zeta Functions},
  pdfauthor={Jayanta Phadikar}
}
\numberwithin{equation}{section}
\allowdisplaybreaks

\theoremstyle{definition}

\theoremstyle{remark}

\theoremstyle{plain}

\DeclareMathOperator{\Li}{Li}

\newcommand{\C}{\mathbb C}
\newcommand{\N}{\mathbb N}
\newcommand{\Nzero}{\mathbb N_0}
\newcommand{\eps}{\varepsilon}

\newcommand{\HH}{\mathcal H}

\newcommand{\Zstrict}{Z^{>}}
\newcommand{\Zstar}{Z^{\star}}
\begin{document}

\begin{frontmatter}

\title{Certified Arbitrary-Precision Evaluation of a Family of Generalized Multiple Zeta Functions}

\author[aff1]{Jayanta Phadikar\texorpdfstring{\corref{cor1}}{}}
\ead{jayantap@wolfram.com}

\address[aff1]{Wolfram Research}

\cortext[cor1]{Corresponding author.}

\begin{abstract}
We describe a certified arbitrary-precision framework for evaluating a family
of generalized multiple zeta functions.  The family includes strict and
weak-star chain sums, ordinary and colored multiple zeta values, affine-base
and polynomial-base variants, and composite levels containing several affine or
polynomial letters with complex coefficients.  The numerical strategy combines
finite-prefix recurrences with two complementary analytic-tail mechanisms:
recursive Euler-Maclaurin expansion of one-variable tails and direct absolute
tail majorants.  The Euler-Maclaurin branch is fast when the relevant suffix
expansions are regular, while the direct-tail branch gives robust certificates
for multi-letter, weak-star, complex-coefficient, and branch-sensitive inputs.
A computation is called certified only when its reported radius is obtained
from a proved analytic bound for the omitted infinite tail.  Strict-disk
colored sums and boundary-color cases with summable absolute majorants are
therefore within the certified scope; conditionally convergent colored cases
whose convergence relies only on non-one unit-modulus oscillation are kept
separate and reported as explicitly non-certified diagnostic outputs unless an
independent analytic remainder bound is available.
\end{abstract}


\end{frontmatter}

\section{Introduction}
\label{sec:introduction}

Throughout the paper, $\N=\{1,2,\ldots\}$.  Exponent parameters may be
complex.  Whenever a base is complex, powers are interpreted after fixing a
branch of the logarithm, and no summand is allowed to contain a zero base at a
positive integer.  Infinite sums are used only under explicit convergence
hypotheses.

Multiple zeta functions are among the central nested sums in experimental and
computational mathematics.  Their classical strict-chain form is
\begin{equation}
  \zeta(s_1,\ldots,s_d)
  =
  \sum_{n_1>\cdots>n_d\ge 1}
  \frac{1}{n_1^{s_1}\cdots n_d^{s_d}} .
  \label{eq:intro-mzv}
\end{equation}
For complex exponents, the ordinary convergence region for
\eqref{eq:intro-mzv} is
\begin{equation}
   \Re(s_1+\cdots+s_k)>k,
   \qquad 1\le k\le d.
   \label{eq:intro-mzv-convergence}
\end{equation}
The case $d=1$ gives the Riemann zeta function.  Higher depth introduces the
strict nesting that is responsible for the stuffle algebra of finite harmonic
sums and for the rich system of relations among multiple zeta values
\citep{Hoffman1992,Zagier1994,BorweinBradleyBroadhurst1997,FlajoletSalvy1998}.

The colored version inserts one decoration at each summation level:
\begin{equation}
  \Li_{s_1,\ldots,s_d}(z_1,\ldots,z_d)
  =
  \sum_{n_1>\cdots>n_d\ge 1}
  \frac{z_1^{n_1}\cdots z_d^{n_d}}
       {n_1^{s_1}\cdots n_d^{s_d}} .
  \label{eq:intro-mpl}
\end{equation}
We regard \eqref{eq:intro-mpl} both as a multiple polylogarithm and as a
colored multiple zeta function.  Root-of-unity colors give alternating and
cyclotomic variants, while general complex colors give the usual
polylogarithmic functions.  For finite strict-chain colored sums we write
\begin{equation}
  \HH_{\mathbf s,\boldsymbol z}(N)
  =
  \sum_{N\ge n_1>\cdots>n_d\ge 1}
  \prod_{j=1}^{d} z_j^{n_j} n_j^{-s_j},
  \qquad
  \Li_{\mathbf s}(\boldsymbol z)=\lim_{N\to\infty}\HH_{\mathbf s,\boldsymbol z}(N)
  \label{eq:intro-finite-colored-harmonic}
\end{equation}
whenever the limit exists.  The uncolored specialization is
$\zeta(\mathbf s)=\Li_{\mathbf s}(1,\ldots,1)$.

The affine-base and polynomial-base zeta families used below, especially the
multi-letter versions in which a single summation level carries several
denominator letters, were introduced in \citet{Phadikar2026MonoidalAlphabets}.
The affine letters have the shape
\begin{equation}
  (\alpha+\beta n)^{-s},
  \label{eq:intro-affine-letter}
\end{equation}
and the polynomial-base letters have the shape
\begin{equation}
  P(n)^{-s},
  \qquad P(n)=c_0+c_1n+\cdots+c_m n^m.
  \label{eq:intro-polynomial-letter}
\end{equation}
The coefficients may be real or complex, subject to the zero-free conditions
needed to define the summand and to keep the asymptotic tail on the chosen
branch.  These classes are natural numerical containers for shifted
Hurwitz-type sums, residue-class decompositions, polynomial zeta functions,
Mathieu-type series, Epstein-Hurwitz-type series, and related sums arising
after partial fraction decomposition or changes of variables
\citep{Elizalde1994,ElizaldeEtAl1994,EieChen1999,Dabrowski2000,PoganyTomovski2006,Matsumoto2006}.
The present paper develops certified arbitrary-precision evaluation methods for
these strict-chain families and introduces their weak-star analogues, where
adjacent summation indices may be equal.

Multiple zeta functions and colored multiple zeta functions occur throughout
the symbolic and numerical theory of special functions.  They arise in Euler
sums, iterated integrals, Mellin transforms, hypergeometric expansions, and
perturbative quantum field theory.  Vermaseren, Bl\"umlein, Kurth, Broadhurst,
Ablinger, Schneider, Phadikar, and others developed symbolic summation tools
and structural frameworks for reducing, transforming, and organizing many of
these objects
\citep{Vermaseren1999,BluemleinKurth1999,BluemleinBroadhurstVermaseren2010,Ablinger2014,Schneider2014,Phadikar2026MonoidalAlphabets}.
In particular, \citet{Phadikar2026MonoidalAlphabets} develops a general sum
formalism that includes colored, affine, and polynomial-base letters.  Those
methods address the symbolic and structural questions: when can a nested sum be
organized, transformed, or reduced to a smaller basis of known functions or
constants?  The present paper is concerned with the complementary numerical
question: given a convergent nested zeta-type sum, how can it be evaluated to
arbitrary precision with reliable control of the infinite tail?

Direct summation is usually inefficient.  In depth $d$, the strict-chain or
weak-star region contains many lattice points, and convergence may be only
polynomial.  High-precision numerical evaluation has therefore been essential
in the subject.  Borwein, Bradley, and Broadhurst used extensive
high-precision computations in their compendium of Euler-Zagier sums
\citep{BorweinBradleyBroadhurst1997}.  Crandall developed fast evaluation
methods for multiple zeta sums \citep{Crandall1998}.  Vollinga and Weinzierl
gave algorithms for numerical evaluation of multiple polylogarithms and
implemented them using arbitrary-precision arithmetic in \textsc{GiNaC}
\citep{VollingaWeinzierl2005}.  For rigorous special-function computation,
Johansson used Euler-Maclaurin summation for high-precision evaluation of the
Hurwitz zeta function and its derivatives, and the \textsc{Arb} library
provided an efficient midpoint-radius framework for propagating numerical
error bounds \citep{JohanssonHurwitzZeta2014,JohanssonArb2017}.

The strategy developed here keeps the nested-sum structure explicit.  A finite
prefix is accumulated by dynamic recurrences, and the omitted infinite part is
then controlled recursively.  In certified regimes, acceptance is based on an
analytic inequality of the form $|Z-V|\le E$, not on agreement between two
numerical truncations alone.  One branch expands the tail by Euler-Maclaurin,
reducing retained terms to Hurwitz tails and to Lerch-type tails in strict-disk
colored cases.  A second branch uses direct absolute majorants of the omitted
region.  The two branches are complementary: Euler-Maclaurin is often much
faster for regular power-like tails, while direct absolute bounds are robust
for multi-letter and branch-sensitive inputs where the reduction to suffix
constants is less favorable.

The term \emph{certified} is used conservatively.  A numerical value is called
certified only when it is accompanied by an analytic radius that bounds the
omitted infinite tail.  Boundary-color calculations that rely only on
unit-modulus oscillatory cancellation are therefore separated from the
certification claim unless an independent analytic remainder bound is
available.

The paper is organized as follows.  Section~\ref{sec:family} defines the
strict and weak-star families, including colored, affine, polynomial-base, and
multi-letter levels.  Section~\ref{sec:tails} gives the finite-prefix
recurrences, the two tail mechanisms, and the scope of the certification
claim.  Section~\ref{sec:implementation} describes the implementation, its
public interface, method selection, and returned certification data.
Section~\ref{sec:conclusion} concludes with limitations and future directions.

\section{The Family of Generalized Sums}
\label{sec:family}

Let $w_j:\N\to\C$ be a one-variable level weight.  The strict affine and
polynomial-base variants in this section are the affine and polynomial-base
multi-letter families discussed in the introduction; the weak-star sums are the
star analogues introduced here.  The strict and weak-star chain sums associated
with $\mathbf w=(w_1,\ldots,w_d)$ are
\begin{equation}
  \Zstrict(\mathbf w)
  =
  \sum_{n_1>\cdots>n_d\ge 1}
  \prod_{j=1}^{d} w_j(n_j),
  \label{eq:strict-general-weight}
\end{equation}
and
\begin{equation}
  \Zstar(\mathbf w)
  =
  \sum_{n_1\ge\cdots\ge n_d\ge 1}
  \prod_{j=1}^{d} w_j(n_j).
  \label{eq:star-general-weight}
\end{equation}
The strict sum is the usual multiple-zeta chain.  The weak-star sum is not
computed by converting it to a linear combination of strict sums; it is a
separate chain geometry with its own prefix recurrence and tail bounds.

The ordinary and colored power-letter cases are obtained from
\begin{equation}
  w_j(n)=x_j^n n^{-s_j}.
  \label{eq:power-level-weight}
\end{equation}
Here $x_j=1$ for ordinary multiple zeta functions, and arbitrary complex
colors $x_j$ give colored multiple zeta functions or multiple polylogarithms.
We use ``strict-disk'' for colors with modulus less than one and
``boundary-color'' for unit-modulus colors.  A boundary color different from
one supplies oscillatory cancellation unless it is neutralized by a cumulative
color product, or by a collapsed block in the weak-star case.

The affine multi-letter case allows a product of affine denominator letters at
a single level:
\begin{equation}
  w_j(n)
  =
  x_j^n
  \prod_{\ell=1}^{m_j}
  (\alpha_{j\ell}+\beta_{j\ell}n)^{-s_{j\ell}}.
  \label{eq:affine-multiletter-weight}
\end{equation}
There is one color $x_j$ per level, not one color per denominator letter.  The
single-letter affine case is $m_j=1$ for every $j$.  Positive real affine data
are included as the special case $\alpha_{j\ell}\ge0$ and $\beta_{j\ell}>0$.
For complex affine data one assumes
\begin{equation}
  \beta_{j\ell}\ne 0,
  \qquad
  \alpha_{j\ell}+\beta_{j\ell}n\ne0
  \quad(n\in\N),
  \label{eq:affine-zero-free}
\end{equation}
so that every summand is defined.

The polynomial-base multi-letter case similarly has
\begin{equation}
  w_j(n)
  =
  x_j^n
  \prod_{\ell=1}^{m_j}P_{j\ell}(n)^{-s_{j\ell}},
  \qquad
  P_{j\ell}(n)=\sum_{r=0}^{M_{j\ell}}c_{j\ell r}n^r.
  \label{eq:polynomial-multiletter-weight}
\end{equation}
The coefficient convention is increasing in powers of $n$.  Positive
polynomial bases, $P_{j\ell}(n)>0$ for $n\in\N$, give a real-denominator
subfamily.  Complex polynomial bases are allowed when the leading coefficient
is nonzero and no positive integer is a root:
\begin{equation}
  c_{j\ell M_{j\ell}}\ne0,
  \qquad
  P_{j\ell}(n)\ne0
  \quad(n\in\N).
  \label{eq:polynomial-zero-free}
\end{equation}
Degree-one polynomial bases are affine bases.

The convergence tests are most compactly stated in terms of effective level
exponents.  For power and affine letters set
\begin{equation}
  \sigma_j=
  \begin{cases}
    s_j, & \text{power letter},\\[2mm]
    \sum_{\ell=1}^{m_j}s_{j\ell}, & \text{affine multi-letter level},
  \end{cases}
  \label{eq:effective-power-affine-exponents}
\end{equation}
and for polynomial-base levels set
\begin{equation}
  \sigma_j=\sum_{\ell=1}^{m_j}M_{j\ell}s_{j\ell}.
  \label{eq:effective-polynomial-exponents}
\end{equation}
These exponents describe the power-law decay of the level weights.  For
uncolored strict chains, the ordinary convergence condition is
\begin{equation}
  \Re(\sigma_1+\cdots+\sigma_k)>k,
  \qquad 1\le k\le d.
  \label{eq:effective-strict-convergence}
\end{equation}
For colored strict chains, let
\begin{equation}
  X_k=x_1x_2\cdots x_k,
  \qquad
  S_k=\sigma_1+\cdots+\sigma_k,
  \qquad
  r_k=\#\{1\le j\le k: X_j=1\}.
  \label{eq:colored-cumulative-data}
\end{equation}
A standard sufficient-and-necessary convergence test for the colored nested
series in this setting is
\begin{equation}
  |X_k|<1
  \quad\text{or}\quad
  \bigl(|X_k|=1\ \text{and}\ \Re(S_k)>r_k\bigr),
  \qquad 1\le k\le d,
  \label{eq:colored-strict-convergence}
\end{equation}
together with exclusion of any prefix satisfying $|X_k|>1$.
The uncolored criterion \eqref{eq:effective-strict-convergence} is recovered
from \eqref{eq:colored-strict-convergence} when all $x_j=1$.

For weak-star chains, one must also consider the strata where adjacent indices
are equal.  Let $\pi=(p_1,\ldots,p_r)$ be a composition of $d$, and collapse
successive blocks of lengths $p_1,\ldots,p_r$.  The collapsed level exponent
is the sum of the $\sigma_j$ in that block, and the collapsed color is the
product of the colors in that block.  The weak-star chain is convergent when
every such contiguous block collapse satisfies the corresponding strict-chain
criterion.  In the uncolored case this means that every collapsed exponent
vector satisfies \eqref{eq:effective-strict-convergence}; in the colored case
it means that every collapsed colored vector satisfies
\eqref{eq:colored-strict-convergence}.

\section{Finite Prefixes and Tail Certification}
\label{sec:tails}

This section describes the numerical core at the mathematical level.  The
common pattern is to compute a finite prefix directly at the chosen working
precision and then control the omitted infinite tail.  The strict and
weak-star prefixes use different dynamic recurrences, while the tails are
handled by either Euler-Maclaurin expansion or direct absolute majorants.

\subsection{Finite strict and weak-star prefixes}
\label{subsec:prefixes}

For a strict chain define
\begin{equation}
  H_i^{>}(N)
  =
  \sum_{N\ge m_i>\cdots>m_d\ge1}
  \prod_{j=i}^{d}w_j(m_j),
  \qquad
  H_{d+1}^{>}(N)=1.
  \label{eq:strict-prefix-suffixes}
\end{equation}
Then
\begin{equation}
  H_i^{>}(N)=H_i^{>}(N-1)+w_i(N)H_{i+1}^{>}(N-1),
  \qquad 1\le i\le d,
  \label{eq:strict-prefix-recurrence}
\end{equation}
with $H_i^{>}(0)=0$ for $i\le d$.  Thus the strict prefix with top index at
most $N$ is $H_1^{>}(N)$.

For a weak-star chain define
\begin{equation}
  H_i^{\star}(N)
  =
  \sum_{N\ge m_i\ge\cdots\ge m_d\ge1}
  \prod_{j=i}^{d}w_j(m_j),
  \qquad
  H_{d+1}^{\star}(N)=1.
  \label{eq:star-prefix-suffixes}
\end{equation}
The update is
\begin{equation}
  H_i^{\star}(N)=H_i^{\star}(N-1)+w_i(N)H_{i+1}^{\star}(N),
  \qquad 1\le i\le d.
  \label{eq:star-prefix-recurrence}
\end{equation}
The occurrence of $H_{i+1}^{\star}(N)$, rather than
$H_{i+1}^{\star}(N-1)$, is precisely the weak-chain correction allowing
$m_i=m_{i+1}$.

\subsection{Euler-Maclaurin tails}
\label{subsec:euler-maclaurin-tails}

The basic one-variable tail is the Hurwitz tail
\begin{equation}
  \sum_{n>N}n^{-a}=\zeta(a,N+1).
\end{equation}
For $a\ne1$ the Euler-Maclaurin expansion gives, for an integer $K\ge1$,
\begin{equation}
  \zeta(a,N+1)
  =
  \frac{(N+1)^{1-a}}{a-1}
  +\frac{(N+1)^{-a}}{2}
  +\sum_{j=1}^{K}
    \frac{B_{2j}}{(2j)!}(a)_{2j-1}(N+1)^{-a-2j+1}
  +\mathcal R_K(a,N),
  \label{eq:hurwitz-em-tail}
\end{equation}
where $B_{2j}$ are Bernoulli numbers.  A useful mixed form is obtained by
multiplying an inner Hurwitz tail by an outer power.  For
\begin{equation}
  S_{p,a}(N)=\sum_{n>N}n^{-p}\zeta(a,n+1),
\end{equation}
one obtains the approximation
\begin{align}
  S_{p,a}(N)
  &\approx
  \frac{1}{a-1}\zeta(p+a-1,N+1)
  +\frac12\zeta(p+a,N+1) \notag\\
  &\quad+
  \sum_{j=1}^{K}
  \frac{B_{2j}}{(2j)!}(a)_{2j-1}
  \zeta(p+a+2j-1,N+1).
  \label{eq:power-hurwitz-tail}
\end{align}
When
\begin{equation}
  \Re(a)+2K-1>0,
  \qquad
  \Re(p)+\Re(a)+2K-1>1,
  \label{eq:power-hurwitz-bound-conditions}
\end{equation}
the discarded part is bounded by
\begin{equation}
  |\mathcal R_K(p,a,N)|
  \le
  \frac{2\zeta(2K)}{(2\pi)^{2K}}
  \frac{|(a)_{2K}|}{\Re(a)+2K-1}
  \zeta(\Re(p)+\Re(a)+2K-1,N+1).
  \label{eq:power-hurwitz-bound}
\end{equation}
This estimate is the basic analytic certificate for ordinary power tails.

The nested recursion may be described in terms of suffix expansions.  A suffix
finite sum has the form
\begin{equation}
  H_{i+1}(n-1)
  =
  C_{i+1}
  +\sum_r c_r n^{-\alpha_r}
  +R_{i+1}(n),
  \label{eq:suffix-expansion}
\end{equation}
where $C_{i+1}$ is the limiting suffix constant and $R_{i+1}(n)$ is bounded
by a finite list of positive power or power-logarithmic majorants.  Substitution
into the next level gives
\begin{equation}
  \sum_{n>N}w_i(n)H_{i+1}(n-1)
  =
  C_{i+1}\sum_{n>N}w_i(n)
  +\sum_r c_r\sum_{n>N}w_i(n)n^{-\alpha_r}
  +\sum_{n>N}w_i(n)R_{i+1}(n).
  \label{eq:tail-recursion-general}
\end{equation}
For power letters, the retained sums in \eqref{eq:tail-recursion-general}
are Hurwitz or Lerch tails.  For affine and polynomial-base letters, the
level weight is first expanded as a finite asymptotic power series and the
remainder is bounded by a positive power majorant.

For an affine letter, the large-$n$ expansion is
\begin{equation}
  (\alpha+\beta n)^{-s}n^{-\gamma}
  =
  \beta^{-s}
  \sum_{r=0}^{R}\binom{-s}{r}
  \left(\frac{\alpha}{\beta}\right)^r
  n^{-s-\gamma-r}
  +\mathcal E_R(n),
  \label{eq:affine-binomial-expansion}
\end{equation}
valid beyond a cutoff at which the chosen logarithm branch is stable.  In a
multi-letter level, the finite expansions of the letters are multiplied
before the one-variable tail is evaluated.  The omitted product terms and the
individual remainders are collected into positive power majorants.

For a polynomial base of degree $m$,
\begin{equation}
  P(n)=c_m n^m\left(1+q_1n^{-1}+\cdots+q_mn^{-m}\right),
\end{equation}
and hence, beyond a branch-safe cutoff,
\begin{equation}
  P(n)^{-s}n^{-\gamma}
  =
  \sum_{r=0}^{R}a_r n^{-ms-\gamma-r}+\mathcal E_R(n).
  \label{eq:polynomial-asymptotic}
\end{equation}
The coefficients $a_r$ come from the Taylor expansion of
$(1+q_1x+\cdots+q_mx^m)^{-s}$ at $x=0$.  A Cauchy estimate on a disk in the
$x=1/n$ plane gives a computable remainder bound of the form
\begin{equation}
  |\mathcal E_R(n)|\le C_R n^{-\Re(ms+\gamma)-R-1},
  \qquad n\ge A,
  \label{eq:polynomial-cauchy-majorant}
\end{equation}
provided the disk avoids zeros of the normalized polynomial factor and remains
within the selected logarithm branch.  Multi-letter polynomial levels are
handled by multiplying the retained level expansions and summing the resulting
power tails.

When an intermediate exponent in the Euler-Maclaurin recursion approaches a
pole, the expression can pass through quantities such as $(a-1)^{-1}$ even
though the original nested sum is convergent.  A finite-part construction
resolves this situation.  The exponents are perturbed in a fixed generic
direction,
\begin{equation}
  (s_1,\ldots,s_d)
  \mapsto
  (s_1,\ldots,s_d)+\eps(2^{-1},2^{-2},\ldots,2^{-d}),
  \label{eq:finite-part-direction}
\end{equation}
the tail expression is expanded as a Laurent series in $\eps$, and the
constant term is retained:
\begin{equation}
  \operatorname{FP}_{\eps=0}F(\eps)=[\eps^0]F(\eps).
  \label{eq:finite-part}
\end{equation}
The same principle applies to composite levels by perturbing at the level-block
scale.  In the certified Euler-Maclaurin branches, the retained finite part is
accompanied by propagated remainder bounds.  Parameter-stability comparisons
are useful diagnostics, but they are not used as substitutes for analytic
remainder bounds in the certification statements.

\subsection{Direct absolute tail majorants}
\label{subsec:direct-tail}

The direct-tail method avoids suffix constants and symbolic reductions.  It
uses positive majorants for the summand and then bounds the omitted chain
region directly.  Suppose each level admits a bound
\begin{equation}
  |w_j(n)|\le C_j q_j^n n^{-\rho_j},
  \qquad 0\le q_j\le1,
  \label{eq:level-majorant}
\end{equation}
possibly after replacing a complex power by the elementary estimate
\begin{equation}
  |z^{-s}|\le e^{\pi|\Im s|}|z|^{-\Re s}
  \label{eq:complex-power-majorant}
\end{equation}
for the principal logarithm.  Affine and polynomial levels satisfy
\eqref{eq:level-majorant} beyond a computable cutoff, with a finite correction
for the initial range.

The suffix of a weak or strict chain with power-type majorants can be bounded
recursively by a finite list of power-logarithmic terms
\begin{equation}
  C n^{\lambda}(\log n)^\ell,
  \qquad C\ge0,
  \quad \ell\in\Nzero.
  \label{eq:power-log-term}
\end{equation}
The first omitted level is then bounded by sums of the form
\begin{equation}
  \sum_{n>N} n^{-\rho+\lambda}(\log n)^\ell.
  \label{eq:power-log-tail}
\end{equation}
When $\rho-\lambda>1$, these tails are controlled by a finite initial
correction plus the integral
\begin{equation}
  \int_N^\infty x^{-\rho+\lambda}(\log x)^\ell\,dx.
  \label{eq:power-log-integral}
\end{equation}
The integral has an explicit finite expression in powers of $\log N$ divided
by powers of $\rho-\lambda-1$.

If a color satisfies $|q|<1$, geometric decay replaces
\eqref{eq:power-log-integral} by a Lerch-type tail.  A basic bound is
\begin{equation}
  \sum_{n>N}|q|^n n^{-\rho}
  =|q|^{N+1}\Phi(|q|,\rho,N+1),
  \label{eq:lerch-tail-bound}
\end{equation}
and logarithmic factors are absorbed by inequalities of the form
$(\log n)^\ell\le (\ell/(e\beta))^\ell n^{\beta}$ with a small positive
$\beta$.  This gives certified absolute bounds in the strict disk.  If
$|q|=1$, the direct-tail certificate deliberately ignores oscillatory
cancellation and therefore succeeds only when the remaining power-logarithmic
absolute majorant is summable.  Thus an alternating example such as
$\sum_{n\ge1}(-1)^n n^{-2}$ can be certified through the absolute majorant
$\sum n^{-2}$, whereas a conditionally convergent boundary example such as
$\sum_{n\ge1}(-1)^n n^{-1/2}$ is outside the direct absolute-tail certified
class.  For mixed color data the geometric and power-logarithmic majorants are
combined level by level; certification never depends on unproved cancellation
between terms.

The direct-tail certificate is usually less asymptotically sharp than a
successful Euler-Maclaurin expansion, but it is structurally robust.  It works
without assigning values to divergent suffix constants, applies naturally to
weak-star chains, and treats multi-letter affine or polynomial levels as a
single composite weight.

\subsection{Scope of certification}
\label{subsec:certified-regimes}

A branch is called certified in this paper only when it produces a value $V_N$
and a positive radius $E_N$ satisfying
\begin{equation}
  |Z-V_N|\le E_N.
  \label{eq:certified-output}
\end{equation}
The radius must come from a proved bound for the omitted analytic tail:
Euler-Maclaurin remainders, affine binomial remainders, polynomial Cauchy
remainders, direct absolute tail majorants, or certified finite-part
combinations of these bounds.  Floating-point roundoff is controlled by guard
precision rather than by a full directed interval-arithmetic proof, so the
certificate asserted here is an analytic truncation certificate for the
mathematical tail.

For colored sums, the distinction is not whether a color lies on the unit
circle, but whether the selected branch proves a tail bound.  Strict-disk data
are certified by absolute Lerch-type majorants.  Boundary-color data with
non-one unit-modulus factors are also certified when oscillation can be ignored
and the absolute majorant is summable, or when another explicit
Euler-Maclaurin/Lerch remainder bound is available.  This includes any case in
which the relevant cumulative color products, and in the weak-star case the
collapsed block products, still lead to a summable certified majorant.

The cases excluded from the certification claim are the conditional boundary
cases in which the remaining tail is known only through cancellation from
non-one unit-modulus factors.  For such inputs, oscillatory Hurwitz/Lerch-type
formulas may be evaluated as numerical extensions: two nearby truncation or
expansion choices are compared, and the observed difference is reported as a
stability indicator.  This diagnostic quantity is not an $E_N$ satisfying
\eqref{eq:certified-output}.

\section{Wolfram Language Implementation}
\label{sec:implementation}

The implementation is supplied as a Wolfram Language package file.  From an
examples notebook in the same directory, it can be loaded by evaluating
\begin{WLExample}
Get[FileNameJoin[{NotebookDirectory[],  "ArbitraryPrecisionZetaFunctions.wl"}]];
\end{WLExample}

The package exposes twelve public numerical evaluators: six strict-chain
functions and six weak-star functions.
\begin{center}
\begin{adjustbox}{max width=\textwidth}
\begin{tabular}{@{}lll@{}}
\toprule
\textbf{Mathematical class} & \textbf{Strict-chain evaluator} & \textbf{Weak-star evaluator} \\
\midrule
ordinary multiple zeta & \texttt{NMultipleZeta} & \texttt{NMultipleZetaS} \\
colored multiple zeta / multiple polylogarithm & \texttt{NMultiplePolyLog} & \texttt{NMultiplePolyLogS} \\
multiple affine zeta & \texttt{NMultipleAffineZeta} & \texttt{NMultipleAffineZetaS} \\
colored multiple affine zeta / affine polylogarithm & \texttt{NMultipleAffinePolyLog} & \texttt{NMultipleAffinePolyLogS} \\
multiple polynomial-base zeta & \texttt{NMultiplePolyBaseZeta} & \texttt{NMultiplePolyBaseZetaS} \\
colored multiple polynomial-base zeta / polynomial-base polylogarithm & \texttt{NMultiplePolyBasePolyLog} & \texttt{NMultiplePolyBasePolyLogS} \\
\bottomrule
\end{tabular}
\end{adjustbox}
\end{center}
The suffix \texttt{S} denotes the weak-star chain
$n_1\ge\cdots\ge n_d\ge1$.  Star sums are accumulated directly; they are not
converted into a sum of strict-chain values.

The input conventions are as follows.
\begin{center}
\begin{adjustbox}{max width=\textwidth}
\begin{tabular}{@{}p{0.24\textwidth}p{0.70\textwidth}@{}}
\toprule
\textbf{Input} & \textbf{Convention} \\
\midrule
exponent list & \texttt{\{s1,...,sd\}} for one letter per level. \\
colors & \texttt{\{x1,...,xd\}}, one color per level.  The all-one color case delegates to the corresponding zeta evaluator. \\
affine bases & \texttt{\{\{a1,b1\},...,\{ad,bd\}\}} represents $(a_j+b_j n)^{-s_j}$ at level $j$. \\
polynomial bases & coefficient lists are increasing-power lists: \texttt{\{c0,c1,...,cm\}} represents $c_0+c_1n+\cdots+c_mn^m$. \\
multi-letter affine or polynomial levels & exponent blocks such as \texttt{\{\{s11,s12\},\{s21\},...\}} are matched with base blocks of the same lengths; all letters in a block are multiplied at that summation level. \\
complex coefficients & affine and polynomial coefficient data may be complex when the no-positive-integer-zero tests pass; powers are taken on the principal branch. \\
\bottomrule
\end{tabular}
\end{adjustbox}
\end{center}
A singleton level should use the ordinary single-letter form; a genuine
multi-letter block is used when a level contains two or more denominator
letters.

All public evaluators accept an optional requested decimal precision; when it
is omitted, the default is 10 digits.  With \texttt{"ReturnCertificationData" -> True}, the result is an
association rather than just the numerical value.  The common keys include
\texttt{"Value"}, \texttt{"ErrorRadius"}, \texttt{"N"}, \texttt{"R"},
\texttt{"K"}, \texttt{"Depth"}, \texttt{"WorkingPrecision"},
\texttt{"Branch"}, and \texttt{"Certification"}.  Depending on the family,
additional metadata may include \texttt{"SelectedMethod"},
\texttt{"MethodSelection"}, \texttt{"LevelLetterCounts"},
\texttt{"Degrees"}, \texttt{"CoefficientDomain"}, and
\texttt{"CoefficientConvention"}.  By default only the value stored under
\texttt{"Value"} is returned.

The public method option has four choices.
\begin{center}
\begin{adjustbox}{max width=\textwidth}
\begin{tabular}{@{}p{0.25\textwidth}p{0.69\textwidth}@{}}
\toprule
\textbf{Method} & \textbf{Meaning} \\
\midrule
\texttt{Automatic} & chooses the branch expected to be fastest among the applicable analytic-tail methods, and usually tries the alternate certification-capable branch if the preferred one cannot prove the requested digits.  A result is treated as certified only when the returned metadata identifies an analytic tail bound.  The exploratory finite-prefix mode is never selected automatically. \\
\texttt{"EulerMaclaurin"} & uses the recursive Euler-Maclaurin tail machinery, including affine binomial tails, polynomial asymptotic tails, and finite-part handling when the required analytic remainder bounds are available.  Oscillatory boundary-color tails are reported as certified only when such bounds are available; otherwise they are marked as explicitly non-certified diagnostics. \\
\texttt{"DirectTail"} & uses a finite prefix plus a certified absolute tail majorant.  This method is especially useful for multi-letter levels, complex coefficient data, strict-disk colored sums, absolutely summable unit-modulus color cases, weak-star sums, and cases where suffix constants are divergent or expensive.  It does not rely on cancellation from unit-modulus colors. \\
\texttt{"CrudeSum"} & returns a fast, machine-precision, 10000-term finite prefix.  It is not certified, reports \texttt{Missing["NotCertified"]} as the error radius when certification data are requested, and is not selected by \texttt{Automatic}. \\
\bottomrule
\end{tabular}
\end{adjustbox}
\end{center}
For \texttt{Automatic} calls, successful data record the chosen branch through
\texttt{"SelectedMethod"} and \texttt{"MethodSelection"}.  The
\texttt{"Certification"} field is part of the contract: only entries described
as certified analytic truncation-error data use \texttt{"ErrorRadius"} as an
analytic tail radius, whereas entries described as ``internal stability check
only'' use \texttt{"ErrorRadius"} as a stability difference and should not be
used as rigorous enclosures.

The principal control options are summarized below.
\begin{center}
\begin{adjustbox}{max width=\textwidth}
\begin{tabular}{@{}lll@{}}
\toprule
\textbf{Option} & \textbf{Typical default} & \textbf{Purpose} \\
\midrule
\texttt{"StartN"} & \texttt{20} or \texttt{Automatic} & initial finite-prefix cutoff. \\
\texttt{"MaxN"} & family dependent & largest cutoff attempted during adaptive doubling. \\
\texttt{"StartR"}, \texttt{"MaxR"} & \texttt{8}, \texttt{48} & affine or polynomial product-expansion order range. \\
\texttt{"StartK"}, \texttt{"MaxK"} & \texttt{4}, \texttt{18} or family dependent & Euler-Maclaurin order range for shifted power tails; accepted by all public evaluators. \\
\texttt{"GuardDigits"} & \texttt{5} & extra working digits used to reduce roundoff sensitivity. \\
\texttt{"SafetyDigits"} & \texttt{2} & extra tolerance margin in certification tests. \\
\texttt{Method} & \texttt{Automatic} & branch selector described above. \\
\texttt{"ReturnCertificationData"} & \texttt{False} & whether to return the full association. \\
\bottomrule
\end{tabular}
\end{adjustbox}
\end{center}
For certified branches, the adaptive search uses the accepted controls $N$, $R$,
and $K$ according to the selected method until the analytic error radius is
below the requested tolerance or
the specified maxima are exhausted.  Failure is a controlled outcome: it means
that the selected branch and parameter bounds did not prove the requested
precision, not necessarily that the mathematical series is divergent.  For
boundary-color diagnostic branches, the reported radius has the separate
non-certified meaning described in Section~\ref{subsec:certified-regimes}.

The ordinary multiple zeta evaluator first checks the convergence region
\eqref{eq:intro-mzv-convergence}.  The Euler-Maclaurin method then chooses
between a non-singular recursive Hurwitz-Euler-Maclaurin branch and a
finite-part branch when intermediate resonances are detected.  The direct-tail
method uses an absolute power-log tail radius.  The weak-star ordinary
evaluator uses the corresponding weak-chain recurrence and weak-tail bounds.

The colored power-letter evaluators detect the all-colors-one case and
redirect it to the appropriate zeta evaluator.  Strict-disk color inputs use
absolute Lerch majorants and are certified.  Boundary inputs with non-one
unit-modulus colors are split according to the proof available for the tail:
when an absolute or analytic majorant proves the remainder, the result is
certified; when convergence depends only on oscillatory cancellation, the
Euler-Maclaurin/Lerch-type branch is returned only as an explicitly
non-certified diagnostic.  In the latter case the reported radius is a
comparison estimate, not an analytic enclosure.

The affine evaluators support positive real data and complex coefficient data.
For complex affine coefficients, the validator requires nonzero slopes and no
positive-integer zero of $a+bn$.  The Euler-Maclaurin path uses branch-safe
binomial tails when every needed suffix is convergent; otherwise the certified
direct-tail fallback can be used.  Multi-letter affine input is handled as a
composite product at each level, with the number of letters per level recorded
in the returned metadata.

The polynomial-base evaluators use increasing-power coefficient lists.  Complex
polynomial coefficients are supported when leading coefficients are nonzero and
there are no positive-integer roots.  Degree-one polynomial inputs are routed
through the affine machinery when this is sharper.  Higher-degree polynomial
Euler-Maclaurin tails use branch-safe asymptotic expansions and Cauchy
remainder bounds; the direct-tail path uses complex polynomial power-log
majorants.  For low requested digit counts in complex polynomial-base cases, an
inexpensive low-order Euler-Maclaurin prepass is tried before the heavier
parameter grid.  Multi-letter polynomial levels are treated as composite level
products rather than expanded into separate one-letter zeta calls.

The supplementary examples notebook exercises the full public interface: strict
and weak-star sums, ordinary and colored power-letter cases, affine and
polynomial-base families, complex coefficients, multi-letter levels, explicit
method choices, certification-data output, and convergence-gate failures.

\clearpage
\subsection{Sample evaluations}
\label{subsec:sample-evaluations}

The following basic calls illustrate the strict-chain evaluators at 10 decimal
digits.
\begin{WLExample}
In[10]:= NMultipleZeta[{2, 3, 4}, 10]
Out[10]= 0.06781184624

In[11]:= NMultiplePolyLog[{3, 2}, {-1/2, 1/4}, 10]
Out[11]= 0.006791728023

In[12]:= NMultipleAffineZeta[{2, 3, 4},
  {{1, 1}, {0, 2}, {2, 3}}, 10]
Out[12]= 0.00001036151260

In[13]:= NMultipleAffinePolyLog[{3 - I/5, 2 + I/7,
  2 + I/8}, {(1 + I)/5, (2 + I)/8, (-1 - I)/7},
  {{1, 1}, {2, 1}, {0, 3}}, 10]
Out[13]= 2.267783230*10^-8 + 3.646884004*10^-8 I

In[14]:= NMultiplePolyBaseZeta[{2, 3},
  {{1, 1, 1}, {1, 2, 1}}, 10]
Out[14]= 0.0005023069889

In[15]:= NMultiplePolyBasePolyLog[{2 + I/4,
  3 - I/6}, {I/3, (-1 + I)/6}, {{1, 0, 1}, {2, 1, 1}}, 10]
Out[15]= 0.00001055661174 - 0.00001264912921 I
\end{WLExample}

\clearpage
The corresponding weak-star evaluators use the same input conventions and add
the suffix \texttt{S} to the function name.
\begin{WLExample}
In[16]:= NMultipleZetaS[{2, 3, 4}, 10]
Out[16]= 1.755290921

In[17]:= NMultiplePolyLogS[{3, 2}, {-1/2, 1/4}, 10]
Out[17]= -0.1177277992

In[18]:= NMultipleAffineZetaS[{2, 3, 4},
  {{1, 1}, {0, 2}, {2, 3}}, 10]
Out[18]= 0.0001445233460

In[19]:= NMultipleAffinePolyLogS[{3 - I/5, 2 + I/7,
  2 + I/8}, {(1 + I)/5, (2 + I)/8, (-1 - I)/7},
  {{1, 1}, {2, 1}, {0, 3}}, 10]
Out[19]= 9.94344375*10^-6 - 0.00002410122361 I

In[20]:= NMultiplePolyBaseZetaS[{2, 3},
  {{1, 1, 1}, {1, 2, 1}}, 10]
Out[20]= 0.002268031210

In[21]:= NMultiplePolyBasePolyLogS[{2 + I/4,
  3 - I/6}, {I/3, (-1 + I)/6}, {{1, 0, 1}, {2, 1, 1}}, 10]
Out[21]= -0.0001935391563 - 0.0002413492341 I
\end{WLExample}

Options can be fixed explicitly.  For example, the following call requests the
Euler-Maclaurin branch and fixes both the product-expansion and shifted-tail
orders.
\begin{WLExample}
In[22]:= NMultipleZeta[{3, 1, 2}, 12, Method -> "EulerMaclaurin",
  "StartR" -> 4, "MaxR" -> 4, "StartK" -> 1, "MaxK" -> 1]
Out[22]= 0.0792213975364
\end{WLExample}

Certification metadata can be requested when one wants the analytic radius and
the selected branch rather than only the numerical value.
\begin{WLExample}
In[23]:= NMultipleAffineZeta[{2, 3, 4}, {{1, 1}, {0, 2}, {2, 3}}, 10,
  "ReturnCertificationData" -> True]

Out[23]= <|"Value" -> 0.00001036151260,
  "ErrorRadius" -> 1.273136583*10^-13, "N" -> 12, "R" -> 8,
  "K" -> 2, "Depth" -> 3, "WorkingPrecision" -> 15,
  "Branch" -> "ordinary convergent positive-affine certified EM",
  "Certification" ->
   "Certified analytic truncation-error radius for Euler-Maclaurin and \
affine binomial tails; internal rounding error is controlled by \
GuardDigits.",
  "SelectedMethod" -> "EulerMaclaurin",
  "MethodSelection" -> "Automatic"|>
\end{WLExample}

\clearpage
\section{Conclusion}
\label{sec:conclusion}

We presented a certified arbitrary-precision framework for a family of
generalized multiple zeta functions.  The family includes ordinary and colored
multiple zeta values, affine and polynomial-base multi-letter strict-chain
families, the weak-star analogues introduced here, and complex coefficient data
under zero-free and branch-safe hypotheses.

The numerical method combines dynamic finite-prefix accumulation with two
complementary analytic-tail strategies.  Euler-Maclaurin tails are fast in
regular situations and reduce retained terms to Hurwitz or Lerch-type
one-variable objects.  Direct absolute tails are often more robust for
branch-sensitive, multi-letter, weak-star, or divergent-suffix inputs.  The
certification claim is confined to branches with analytic tail radii;
boundary-color diagnostic outputs are deliberately separated from certified
output.

The approach is intended as a numerical complement to symbolic reduction
methods.  It supports testing identities, validating conjectural reductions,
performing integer-relation searches, and exploring constants for which no
closed reduction is known.  Future work includes certified cancellation-aware
bounds for conditionally convergent unit-modulus colored tails, stronger
directed roundoff control, and tighter integration with stuffle reduction and
Shintani-type decompositions.

\paragraph{AI-assisted preparation disclosure}
AI-assisted tools were used in a limited way to improve the clarity and presentation
of selected parts of the manuscript. The mathematical ideas, definitions, results,
proofs, computations, and overall substance are the author's own. All AI-assisted
suggestions were carefully reviewed, and the author remains fully responsible for
the accuracy, originality, and final form of the article.

\section*{Declaration of competing interest}
The author declares no competing interests.


\begin{thebibliography}{99}

\bibitem[Ablinger(2014)]{Ablinger2014}
Ablinger, J., 2014. Computer algebra algorithms for special functions in
particle physics. Ph.D. thesis, Johannes Kepler University Linz.

\bibitem[Bl\"umlein and Kurth(1999)]{BluemleinKurth1999}
Bl\"umlein, J., Kurth, S., 1999. Harmonic sums and Mellin transforms up to
two-loop order. \emph{Phys. Rev. D} 60, 014018.

\bibitem[Bl\"umlein et~al.(2010)]{BluemleinBroadhurstVermaseren2010}
Bl\"umlein, J., Broadhurst, D.J., Vermaseren, J.A.M., 2010. The multiple zeta
value data mine. \emph{Comput. Phys. Commun.} 181, 582--625.

\bibitem[Borwein et~al.(1997)]{BorweinBradleyBroadhurst1997}
Borwein, J.M., Bradley, D.M., Broadhurst, D.J., 1997. Evaluations of
$k$-fold Euler/Zagier sums: a compendium of results for arbitrary $k$.
\emph{Electron. J. Combin.} 4(2), R5.

\bibitem[Crandall(1998)]{Crandall1998}
Crandall, R.E., 1998. Fast evaluation of multiple zeta sums.
\emph{Mathematics of Computation} 67, 1163--1172.

\bibitem[Dabrowski(2000)]{Dabrowski2000}
Dabrowski, A., 2000. On zeta functions associated with polynomials.
\emph{Bulletin of the Australian Mathematical Society} 61, 455--462.

\bibitem[Duhr et~al.(2012)]{DuhrGanglRhodes2012}
Duhr, C., Gangl, H., Rhodes, J.R., 2012. From polygons and symbols to
polylogarithmic functions. \emph{JHEP} 10, 075.

\bibitem[Eie and Chen(1999)]{EieChen1999}
Eie, M., Chen, K.-W., 1999. A theorem on zeta-functions associated with
polynomials. \emph{Transactions of the American Mathematical Society} 351,
3217--3228.

\bibitem[Elizalde(1994)]{Elizalde1994}
Elizalde, E., 1994. Analysis of an inhomogeneous generalized
Epstein-Hurwitz zeta function with physical applications. \emph{Journal of
Mathematical Physics} 35, 6100--6122.

\bibitem[Elizalde et~al.(1994)]{ElizaldeEtAl1994}
Elizalde, E., Odintsov, S.D., Romeo, A., Bytsenko, A.A., Zerbini, S., 1994.
\emph{Zeta Regularization Techniques with Applications}. World Scientific,
Singapore.

\bibitem[Flajolet and Salvy(1998)]{FlajoletSalvy1998}
Flajolet, P., Salvy, B., 1998. Euler sums and contour integral
representations. \emph{Experiment. Math.} 7, 15--35.

\bibitem[Hoffman(1992)]{Hoffman1992}
Hoffman, M.E., 1992. Multiple harmonic series. \emph{Pacific J. Math.} 152,
275--290.

\bibitem[Johansson(2014)]{JohanssonHurwitzZeta2014}
Johansson, F., 2014. Rigorous high-precision computation of the Hurwitz zeta
function and its derivatives. \emph{Numerical Algorithms} 69, 253--270.

\bibitem[Johansson(2017)]{JohanssonArb2017}
Johansson, F., 2017. Arb: efficient arbitrary-precision midpoint-radius
interval arithmetic. \emph{IEEE Transactions on Computers} 66, 1281--1292.

\bibitem[Lewin(1981)]{Lewin1981}
Lewin, L., 1981. \emph{Polylogarithms and Associated Functions}. North-Holland,
New York.

\bibitem[Matsumoto(2006)]{Matsumoto2006}
Matsumoto, K., 2006. Asymptotic expansions of double zeta-functions of Barnes,
of Shintani, and Eisenstein series. \emph{Nagoya Mathematical Journal} 172,
59--102.

\bibitem[Phadikar(2026)]{Phadikar2026MonoidalAlphabets}
Phadikar, J., 2026. Monoidal alphabets for generalized harmonic sums.
arXiv:2605.21525.

\bibitem[Pogany and Tomovski(2006)]{PoganyTomovski2006}
Pogany, T.K., Tomovski, Z., 2006. On Mathieu-type series. \emph{Integral
Transforms and Special Functions} 17, 881--895.

\bibitem[Remiddi and Vermaseren(2000)]{RemiddiVermaseren2000}
Remiddi, E., Vermaseren, J.A.M., 2000. Harmonic polylogarithms. \emph{Int. J.
Mod. Phys. A} 15, 725--754.

\bibitem[Schneider(2014)]{Schneider2014}
Schneider, C., 2014. Modern summation methods for loop integrals in quantum
field theory: the packages \textsc{Sigma}, \textsc{EvaluateMultiSums} and
\textsc{SumProduction}. \emph{J. Phys. Conf. Ser.} 523, 012037.

\bibitem[Vermaseren(1999)]{Vermaseren1999}
Vermaseren, J.A.M., 1999. Harmonic sums, Mellin transforms and integrals.
\emph{Int. J. Mod. Phys. A} 14, 2037--2076.

\bibitem[Vollinga and Weinzierl(2005)]{VollingaWeinzierl2005}
Vollinga, J., Weinzierl, S., 2005. Numerical evaluation of multiple
polylogarithms. \emph{Comput. Phys. Commun.} 167, 177--194.

\bibitem[Zagier(1994)]{Zagier1994}
Zagier, D., 1994. Values of zeta functions and their applications. In:
Joseph, A., Mignot, F., Murat, F., Prum, B., Rentschler, R. (Eds.),
\emph{First European Congress of Mathematics}, Vol. II. Birkh\"auser, Basel,
pp. 497--512.

\end{thebibliography}
\end{document}